\documentclass[12pt,thmsa]{article}
\usepackage{amsmath}
\usepackage{graphicx}
\usepackage{amsfonts}
\usepackage{amssymb}
\usepackage{mathrsfs}

\newtheorem{thm}{Theorem}
\newtheorem{lem}{Lemma}
\newtheorem{prop}[thm]{Proposition}

\begin{document}

\begin{center}
{\Large \textbf{Asymptotic normality of  a  generalized maximum mean discrepancy estimator}}

\bigskip

Armando Sosthene Kali  BALOGOUN\textsuperscript{a} , Guy Martial  NKIET\textsuperscript{b}  and Carlos OGOUYANDJOU\textsuperscript{a}

\bigskip

\textsuperscript{a}Institut de Math\'ematiques et de Sciences Physiques, Porto Novo, B\'enin.
\textsuperscript{b}Universit\'{e} des Sciences et Techniques de Masuku,  Franceville, Gabon.

\bigskip

E-mail adresses : sosthene.balogoun@imsp-uac.org,  guymartial.nkiet@mathsinfo.univ-masuku.com,  ogouyandjou@imsp-uac.org.

\bigskip
\end{center}

\noindent\textbf{Abstract.} In this paper, we propose an estimator of the  generalized maximum mean discrepancy between several distributions, constructed by modifying a naive estimator. Asymptotic normality is obtained for this estimator both under equality of these distributions  and under the alternative hypothesis.

\bigskip

\noindent\textbf{AMS 1991 subject classifications: }62E20, 46E22.

\noindent\textbf{Key words:} Asymptotic normality;  Kernel method;  Generalized maximum mean discrepancy.
\section{Introduction}
\label{Intro}
 \noindent When adressing the problem of testing whether two distributions are equal on the basis of samples drawn from each of them, Gretton et al. (2007, 2012) introduced the Maximum Mean Discrepancy (MMD) in reproducing kernel Hilbert space. The MMD is used as test statisctic but its asymptotic null distribution is an infinite sum of distributions, and as such it is not easy to use for achieving the testing procedure. For overcoming such drawback, Makigusa and Naito (2020) adopted an approach proposed in Ahmad (1993) consisting in making an  appropriate modification on the test statistic in order to yield asymptotic normality both under the  null hypothesis and under the alternative. However, they only dealt with the problem of testing whether an unknown distribution is equal to a specified one.  So, it may be   of interest to extend their approach to testing for the equality of  two or more unknown distributions.  Recently, Balogoun et al. (2018) introduced the generalized maximum mean discrepancy (GMMD)   in reproducing kernel Hilbert space, that allows one to deal with more than two distributions, and to test wheher these unknown distributions are equal. In this paper, we propose an estimator of the GMMD  constructed by modifying a naive estimator, and we obtain asymptotic normality  for this estimator both under equality of these distributions  and under the alternative hypothesis. The GMMD is recalled in Section 2, and Section 3 is devoted to its estimation and to the main results. All the proofs are postponed in Section 4.

\section{The generalized maximum mean discrepancy}\label{sec2}

\noindent  Let us  consider a reproducing kernel  Hilbert space (RKHS)  $\mathcal{H}$ of functions from a metric space $\mathcal{X}$ to $\mathbb{R}$.  Throughout this paper, we assume that $K$ satisfies the following assumption:

\bigskip

\noindent $(\mathscr{A}_1):$  $\Vert  K\Vert_{\infty}:= \sup\limits_{(x,y)\in \mathcal{X}^2} K(x,y)<+ \infty $.

\bigskip

\noindent For $\ell\in\{1,\cdots,s\}$ with $s\geq 2$, let   $X_\ell$ be a random variable with values into $\mathcal{X}$  and distribution denoted by $\mathbb{P}_\ell$. From ($\mathscr{A}_1$),    $\mathbb{E}(\sqrt{K(X_\ell,X_\ell)})<+\infty$,  hence the kernel mean embeding $m_\ell$ of $\mathbb{P}_\ell$  exists; it is defined by $m_\ell=\mathbb{E}\left(K(X_\ell,\cdot)\right)$. For the case of $s=2$, Gretton et al (2007, 2012) defined  the maximum mean discrepancy (MMD) as the distance between $\mathbb{P}_1$ and $\mathbb{P}_2$ given by:
\begin{equation*}
\textrm{MMD}(\mathbb{P}_1,\mathbb{P}_2):=\Vert m_1-m_2\Vert_\mathcal{H},
\end{equation*}
where $\Vert \cdot\Vert_\mathcal{H}$ denotes the norm induced by the inner product $<\cdot,\cdot>_\mathcal{H}$ of  $\mathcal{H}$. A generalisation of this notion, that allows one to deal with the case of $s>2$, was given in Balogoun et al. (2018) and is recalled below. 

\bigskip

\noindent\textbf{Definition 1.}\textsl{
The generalized maximum mean discrepancy (GMMD) of the distributions  $\mathbb{P}_1, \cdots,\mathbb{P}_s$, related to  and $\eta=\left(\eta_1,\cdots,\eta_s\right)\in ]0,1[^s$  with $\sum_{\ell=1}^s\eta_\ell=1$, is:}
\begin{equation*}
 \textrm{GMMD}^2(\mathbb{P}_1, \cdots,\mathbb{P}_s;\eta)=\sum_{j=1}^{s}\sum_{\underset{\ell\neq j}{\ell=1}}^{s}\eta_\ell\,\textrm{MMD}^2(\mathbb{P}_j,\mathbb{P}_\ell)=\sum_{j=1}^{s}\sum_{\underset{\ell\neq j}{\ell=1}}^{s}\eta_\ell \parallel m_j-m_\ell\parallel_{\mathcal{H}}^2.
\end{equation*}
 
\bigskip

\noindent This  definition recovers that of   MMD that appears to be a particular case  obtained for  $s=2$.  The hypothesis $\mathscr{H}_0\,:\,\mathbb{P}_1= \cdots =\mathbb{P}_s$ can be characterized by means of the GMMD. Indeed, it is easy to check that   this hypothesis  is true if, and only if,   $\textrm{GMMD}(\mathbb{P}_1, \cdots,\mathbb{P}_s;\eta)=0$ for any  $\eta \in ]0,1[^s$.

\section{Estimation of GMMD and asymptotic normality} 
\noindent For any $j\in\{1,\cdots,s\}$, let $X_1^{(j)},\cdots,X_{n_j}^{(j)}\in\mathcal{X}$ be an i.i.d. sample drawn from $\mathbb{P}_j$. We assume that these samples are independent, i.e. $X_i^{(j)}\perp  X_p^{(\ell)}$ for $j\neq\ell$ and any $(i,p)\in\{1,\cdots,n_j\}\times \{1,\cdots,n_\ell\}$, where $\perp$ denotes stochastic independence.
Putting $n=\sum_{j=1}^sn_j$ and $\pi_j=\frac{n_j}{n}$, we make  the folowing assumption:

\bigskip

\noindent $(\mathscr{A}_2):$ For $j\in \{1,\cdots,s\}$,  there exists $\rho_j\in]0,1[$ such that

 $\lim\limits_{n_j\rightarrow +\infty}\bigg\{\sqrt{n}\left(\pi_j- \rho_j\right)\bigg\}=0$.  

\bigskip

\noindent This assumption implies that $\lim\limits_{n_j\rightarrow +\infty} \left(\pi_j\right)=\rho_j$ and $\sum_{j=1}^s \rho_j=1$. Note that it is always possible to take the previous samples so that  $(\mathscr{A}_2)$ holds. Indeed, for any $(\rho_1,\cdots,\rho_s)\in ]0,1[^s$ satisfying $\sum_{j=1}^s \rho_j=1$ and any $n\in\mathbb{N}^\ast$, it suffices to put $n_j=[n\rho_j]$ for $j\in\{1,\cdots,s-1\}$, where $[a]$ denotes the integer part of $a$, and $n_s=n-\sum_{j=1}^{s-1}n_j$.

\bigskip

\noindent Based on the previous samples,  a naive consistent estimator $\widehat{\mathscr{T}}_n$ of the parameter $\mathscr{T}=\textrm{GMMD}^2(\mathbb{P}_1, \cdots,\mathbb{P}_s;\rho)$ (with $\rho=(\rho_1,\cdots,\rho_s)$) is obtained by replacing each $m_j$ by $\widehat{m}_j=n_j^{-1}\sum_{i=1}^{n_j}K(X_i^{(j)},\cdot)$ and $\rho_j$ by $\pi_j$, i.e.
\begin{eqnarray}\label{tn}
\widehat{\mathscr{T}}_n&=&\sum_{j=1}^{s}\sum_{\underset{\ell\neq j}{\ell=1}}^{s}\pi_\ell \parallel \widehat{m}_j-\widehat{m}_\ell\parallel_{\mathcal{H}}^2\nonumber\\ 
&=&\sum_{j=1}^{s}\sum_{\underset{\ell\neq j}{\ell=1}}^{s}\pi_\ell \bigg\{\parallel \widehat{m}_j \parallel_{\mathcal{H}}^2+\parallel  \widehat{m}_\ell\parallel_{\mathcal{H}}^2
-\frac{2}{n_j}\sum_{i=1}^{n_j}<K(X_i^{(j)},\cdot),\widehat{m}_\ell>_{\mathcal{H}}\bigg\}.
\end{eqnarray}
But, although asymptotic normality can be obtained for this estimator, we found that, under $\mathscr{H}_0$, the  asymptotic variance equals $0$, so this statistic cannot be used for testing for equality of the distributions. That is why, following an approach used in Ahmad (1993) and Makigusa and Naito (2020), we propose an estimator  $\widehat{\mathscr{T}}_{n,\gamma}$ obtained by applying weights $k_{i,n_j}(\gamma)$ to the cross-product terms of (\ref{tn}), i.e.
\begin{eqnarray*}
\widehat{\mathscr{T}}_{n,\gamma}&=&\sum_{j=1}^{s}\sum_{\underset{\ell\neq j}{\ell=1}}^{s}\pi_\ell \bigg\{\parallel \widehat{m}_j \parallel_{\mathcal{H}}^2+\parallel  \widehat{m}_\ell\parallel_{\mathcal{H}}^2
-\frac{2}{n_j}\sum_{i=1}^{n_j}k_{i,n_j}(\gamma)<K(X_i^{(j)},\cdot),\widehat{m}_\ell>_{\mathcal{H}}\bigg\}\\ 
&=&\sum_{j=1}^{s}\sum_{\underset{\ell\neq j}{\ell=1}}^{s}\pi_\ell \bigg\{\frac{1}{n_j^2}\sum_{i,p=1}^{n_j}K(X_i^{(j)},X_p^{(j)})+\frac{1}{n_\ell^2}\sum_{i,p=1}^{n_\ell}K(X_i^{(\ell)},X_p^{(\ell)}) \\
&&-\frac{2}{n_jn_\ell}\sum_{i=1}^{n_j}\sum_{p=1}^{n_\ell}k_{i,n_j}(\gamma)\,K(X_i^{(j)},X_p^{(\ell)})\bigg\}.
\end{eqnarray*}
As in Makigus and Naito (2020), the weights $\left(k_{i,r}(\gamma)\right)_{1\leq i\leq r}$ are positive real numbers depending on a parameter $\gamma\in]0,1]$ and satisfying the following assumptions:

\bigskip
\noindent $(\mathscr{A}_3):$ There exists  a strictly positive real number $ \tau $  and an  integer $ n_0 $ such that for all $r> n_0$:
$$r\left\vert\frac{1}{r}\sum_{i=1}^{n_j}k_{i,r}(\gamma)-1\right\vert\leq \tau.$$\\

\noindent $(\mathscr{A}_4):$ There exists  $ c_k $ such that $\max\limits_{1\leq i \leq r}k_{i,r}(\gamma)< c_k$ for all  $r\in\mathbb{N}^\ast$ and $ \gamma \in]0, 1]$.\\

\noindent$(\mathscr{A}_5):$ for any   $ \gamma \in]0, 1]$,
$\lim\limits_{r\rightarrow +\infty}\frac{1}{r}\sum_{i=1}^{r}k^2_{i,r}(\gamma)=k^2(\gamma)> 1$.

\bigskip

\noindent A typical example is given by $k_{i,r}(\gamma)=1+(-1)^i\,\gamma$ (see Ahmad (1993)). Now, we are able to give asymptotic normality for this estimator. Putting $m=\sum_{j=1}^{s}\rho_j\,m_j$ and  $\mu=\sum_{j=1}^{s}m_j $, and considering the functions $\mathcal{U}_j$ and $\mathcal{V}_j$ from $\mathcal{X}$ to $\mathbb{R}$  defined by 
$\mathcal{U}_{j}(x)=<K(x,\cdot)-m_j, (1-2\rho_j +s\rho_j)m_j+\rho_j(m_j-\mu)>_{\mathcal{H}}$ and 
$\mathcal{V}_{j}(x)=<K(x,\cdot)-m_j,m-\rho_jm_j>_{\mathcal{H}}$,   we have:

\begin{thm}\label{t1}
	Assume  that  $(\mathscr{A}_1)$ to $(\mathscr{A}_5)$ hold. Then as $\min\limits_{1\leq j\leq s}(n_j)\rightarrow +\infty$, we have $\sqrt{n}\{\widehat{\mathscr{T}}_{n,\gamma}-\mathscr{T}\}	\stackrel{\mathscr{D}}{\rightarrow} \mathcal{N}\left(0,\sigma^2_{\gamma}\right)$, 
	where $\stackrel{\mathscr{D}}{\rightarrow}$ denotes convergence in distribution,  and $\sigma^2_{\gamma}=\sum_{j=1}^{s}4\rho_j^{-1}\sigma_j^2(\gamma)$ with:
	\begin{eqnarray}\label{sigmaj}
	\sigma_j^2(\gamma)&=&Var\left(\mathcal{U}_{j}(X_1^{(j)})\right) +k^2(\gamma)\,Var\left(\mathcal{V}_{j}(X_1^{(j)})\right) \nonumber\\
&&-2Cov\bigg(\mathcal{U}_{j}(X_1^{(j)}) ,\mathcal{V}_{j}(X_1^{(j)}) \bigg).
	\end{eqnarray}

\end{thm}
\noindent\textbf{Remark 1.} When $ \mathbb{P}_1= \cdots=\mathbb{P}_s$, we have  $m_1=m_2=\cdots=m_k=m$. Thus $\mathcal{U}_{j}(x)=\mathcal{V}_{j}(x)=(1-\rho_j)<K(x,\cdot)-m ,m>_{\mathcal{H}}$, and
	
	\begin{equation*}
	\sigma^2_{\gamma}=4\left(k^2(\gamma)-1\right)\nu^2\sum_{j=1}^{s}4\rho_j^{-1}\left(1-\rho_j\right)^2,
	\end{equation*}
where $\nu^2=Var\left(<K(X_1^{(1)},\cdot)-m ,m>_{\mathcal{H}}\right)=Var\left(<K(X_1^{(1)},\cdot) ,m>_{\mathcal{H}}\right)$. This shows that  $\widehat{\mathscr{T}}_{n,\gamma}$ has asymptotic normality both under $\mathscr{H}_0$ and under the  alternative hypothesis and, cosequently, that it can be used as a test statistic for testing for $\mathscr{H}_0$.

\bigskip

\noindent In the case of $\mathbb{P}_1=\cdots=\mathbb{P}_s$, we can obtain a consistent estimator of  $\sigma^2_{\gamma}$. Indeed, putting $\widehat{m}=\sum_{j=1}^s\pi_j\widehat{m}_j$,
\[
\widehat{\nu}^2_j=\frac{1}{n_j}\sum_{i=1}^{n_j}<K(X_i^{(j)},\cdot) ,\widehat{m}>_{\mathcal{H}}^2-\bigg(\frac{1}{n_j}\sum_{i=1}^{n_j}<K(X_i^{(j)},\cdot) ,\widehat{m}>_{\mathcal{H}}\bigg)^2
\]
and $\widehat{\nu}^2=\sum_{j=1}^s\pi_j\widehat{\nu}_j^2$, we have:

\begin{prop}\label{pro}
	Assume  that $(\mathscr{A}_1)$ and  $(\mathscr{A}_2)$ hold. Then,  as $\min\limits_{1\leq j\leq s}(n_j)\rightarrow +\infty$, the estimator 
$
\widehat{\sigma}^2_{\gamma}=4\left(k^2(\gamma)-1\right)\widehat{\nu}^2\sum_{j=1}^{s}4\pi_j^{-1}\left(1-\pi_j\right)^2
$
is consistent for $\sigma^2_{\gamma}$ under $\mathbb{P}_1=\cdots=\mathbb{P}_s$.
\end{prop}
\section{Proofs}
\subsection{Preliminary result}
\noindent Putting
\begin{eqnarray}
 &&An=\sqrt{n}\sum_{j=1}^{s}\sum_{\underset{\ell\neq j}{\ell=1}}^{s}\,(\pi_\ell-\rho_\ell)\,\Vert \widehat{m}_j- \widehat{m}_\ell\Vert_\mathcal{H},\label{defan}\\
 && B_n=\sqrt{n}\sum_{j=1}^{s}\sum_{\underset{\ell\neq j}{\ell=1}}^{s}\,(\pi_\ell-\rho_\ell)\,\left\{\frac{1}{n_j}\sum_{i=1}^{n_j}\left(k_{i,n_j}(\gamma)-1\right)<K(X_i^{(j)},\cdot),\widehat{m}_\ell>_{\mathcal{H}}\right\},\label{defbn}\\
 && C_n=\sqrt{n}\sum_{j=1}^{s}\sum_{\underset{\ell\neq j}{\ell=1}}^{s}\,\rho_\ell\,<\widehat{m}_j-m_j,\widehat{m}_\ell-m_\ell>_{\mathcal{H}},\label{defcn}\\
 && D_n=\sqrt{n}\sum_{j=1}^{s}\sum_{\underset{\ell\neq j}{\ell=1}}^{s}\,\,\frac{\rho_\ell}{n_j}\sum_{i=1}^{n_j}\left(k_{i,n_j}(\gamma)-1\right)\left\{<K(X_i^{(j)},\cdot),\widehat{m}_\ell-m_\ell>_{\mathcal{H}}+<m_j,m_\ell>_{\mathcal{H}}\right\}\label{defdn},
 \end{eqnarray} 
we have:
\begin{lem}\label{lem1}
	Assume that $(\mathscr{A}_1)$ to  $(\mathscr{A}_5)$ hold. Then $A_n$, $B_n$, $C_n$ and $D_n$ converge in probability to $0$ as $\min\limits_{1\leq j\leq k}(n_j)\rightarrow +\infty$.	
\end{lem}
\noindent\textit{Proof.} First,  for any $(j,\ell)\in\{1,\cdots,k\}^2$, we have
	$
	\left\|\widehat{m}_j-\widehat{m}_\ell\right\|_{\mathcal{H}}\leq\left\|\widehat{m}_j-m_j\right\|_{\mathcal{H}}+\left\|m_j-m_\ell\right\|_{\mathcal{H}}+\left\|m_\ell-\widehat{m}_\ell\right\|_{\mathcal{H}} 
	$
	and   
\begin{equation}\label{diffmoy}
\left\|\widehat{m}_j-m_j\right\|_{\mathcal{H}}=O_P(n_j^{-1/2}). 
\end{equation}
	Since, from assumption  ($\mathscr{A}_2$), $\lim\limits_{n_j\rightarrow +\infty}\sqrt{n}(\pi_j-\rho_j)=0$  and since $n_j^{-1}\rightarrow 0$, we deduce that 
	$\sqrt{n}(\pi_\ell-\rho_\ell)\Vert\widehat{m}_j-\widehat{m}_\ell\Vert_\mathcal{H}^2=o_P(1)$  and, therefore,  $A_n=o_P(1)$. Secondly, puting
	\begin{eqnarray*}
	B_{j,\ell,n}=\sqrt{n}(\pi_\ell-\rho_\ell)\,\left\{\frac{1}{n_j}\sum_{i=1}^{n_j}\left(k_{i,n_j}(\gamma)-1\right)<K(X_i^{(j)},\cdot),\widehat{m}_\ell>_{\mathcal{H}}\right\},\nonumber
	\end{eqnarray*}
	we obtain by using the Cauchy-Schwartz inequality:
	\begin{eqnarray*}\label{bjln}
	\left|B_{j,\ell,n}\right|&\leq&\left|\pi_\ell-\rho_\ell\right|\,\left|\frac{\sqrt{n}}{n_j}\sum_{i=1}^{n_j}\left(k_{i,n_j}(\gamma)-1\right)\right|\,\left\|K(X_i^{(j)},\cdot)\right\|_{\mathcal{H}}\,\left(\left\|\widehat{m}_\ell-m_\ell\right\|_{\mathcal{H}}+\left\|m_\ell\right\|_{\mathcal{H}}\right).
	\end{eqnarray*}
	On the  one hand,  
\begin{equation}\label{kinf}
\left\|K(X_i^{(j)},\cdot)\right\|_{\mathcal{H}}=\sqrt{<K(X_i^{(j)},\cdot),K(X_i^{(j)},\cdot)>_{\mathcal{H}}}=\sqrt{K(X_i^{(j)},X_i^{(j)})}\leq\left\|K\right\|^{1/2}_{\infty}
\end{equation}
and, on the other hand, using the assumption $(\mathscr{A}_3)$, we have
	\begin{eqnarray}
	\left|\frac{\sqrt{n}}{n_j}\sum_{i=1}^{n_j}\left(k_{i,n_j}(\gamma)-1\right)\right|&=&\frac{\pi_j^{-1}}{\sqrt{n}}\left\{n_j \left|\frac{1}{n_j}\sum_{i=1}^{n_j}k_{i,n_j}(\gamma)-1\right|\right\}\leq \frac{\pi_j^{-1}}{\sqrt{n}}\tau .\nonumber 
	\end{eqnarray}	  
	Then, since   $\lim\limits_{n_\ell\rightarrow +\infty}(\pi_\ell-\rho_\ell)=0$ and   $\lim\limits_{n_j\rightarrow +\infty}\pi_j^{-1}=\rho^{-1}_j$, we deduce   from (\ref{diffmoy}) and  the preceding inequalities  that  $B_{j,\ell,n}=o_P(1)$. Hence, from  the equality $ B_n=\sum_{j=1}^{k}\sum_{\underset{\ell\neq j}{\ell=1}}^{k}B_{j,\ell,n}$, we deduce  that  $B_n=o_P(1)$. Thirdly, using Cauchy-Schwartz inequality, we obtain:
	\begin{eqnarray}\label{cn}
	\left|C_n\right|&\leq&\sum_{j=1}^{s}\sum_{\underset{\ell\neq j}{\ell=1}}^{s}\,\rho_\ell\,\sqrt{n}\left\|\widehat{m}_j-m_j\right\|_{\mathcal{H}}\left\|\widehat{m}_\ell-m_\ell\right\|_{\mathcal{H}};
	\end{eqnarray}
From (\ref{diffmoy}),  it follows that
	$\sqrt{n}\left\|\widehat{m}_j-m_j\right\|_{\mathcal{H}}\,\left\|\widehat{m}_\ell-m_\ell\right\|_{\mathcal{H}}=O_P(\frac{\sqrt{n}}{\sqrt{n_j}\sqrt{n_\ell}}),$ 
and since
 $$\lim\limits_{n_j,n_\ell\rightarrow +\infty}\frac{\sqrt{n}}{\sqrt{n_j}\sqrt{n_\ell}}=\lim\limits_{n_j,n_\ell\rightarrow +\infty}\frac{\sqrt{n}}{\sqrt{n_j}}\, \frac{1}{\sqrt{n_\ell}}=\rho_j^{-1/2}\lim\limits_{n_\ell\rightarrow +\infty} \frac{1}{\sqrt{n_\ell}}=0,$$
	we obtain: $
	\sqrt{n}\left\|\widehat{m}_j-m_j\right\|_{\mathcal{H}}\,\left\|\widehat{m}_\ell-m_\ell\right\|_{\mathcal{H}}=o_P(1)
	$. Then,   (\ref{cn}) allows us  to conclude that  $C_n=o_P(1)$. Fouth, using Cauchy-Schwartz inequality and assumption  $(\mathscr{A}_3)$, we obtain
	\begin{eqnarray}
	\left|D_n\right|
&\leq&\sum_{j=1}^{s}\sum_{\underset{\ell\neq j}{\ell=1}}^{s}\rho_\ell\frac{\pi_j^{-1}}{\sqrt{n}}\left\{n_j \left|\frac{1}{n_j}\sum_{i=1}^{n_j}k_{i,n_j}(\gamma)-1\right|\right\}\left\{\left\|K(X_i^{(j)},\cdot)\right\|_{\mathcal{H}}\,\left\|\widehat{m}_\ell-m_\ell\right\|_{\mathcal{H}}+\left\|m_\ell\right\|_{\mathcal{H}}\, \left\|m_j\right\|_{\mathcal{H}}\right\} \nonumber\\
&\leq&\sum_{j=1}^{s}\sum_{\underset{\ell\neq j}{\ell=1}}^{s}\frac{\pi_j^{-1}}{\sqrt{n}}\tau\rho_\ell \left\{\left\|K(X_i^{(j)},\cdot)\right\|_{\mathcal{H}}\,\left\|\widehat{m}_\ell-m_\ell\right\|_{\mathcal{H}}+\left\|m_\ell\right\|_{\mathcal{H}}\, \left\|m_j\right\|_{\mathcal{H}}\right\}. \nonumber
	\end{eqnarray}
This inequality, together with   (\ref{diffmoy}), (\ref{kinf}) and the fact that 
	$\lim\limits_{n_j\rightarrow +\infty}\pi_j^{-1}=\rho^{-1}_j$, allows us  to conclude that  $D_n=o_P(1)$.
\subsection{Proof of Theorem \ref{t1}}
\noindent Clearly, $\sqrt{n}\left(\widehat{\mathscr{T}}_{n,\gamma}-\mathscr{T}\right)=\delta_n+U_n$, where
$
\delta_n=\sqrt{n}\sum_{j=1}^{s}\sum_{\underset{\ell\neq j}{\ell=1}}^{s}(\pi_\ell-\rho_\ell)\,\widehat{\Gamma}^{(n)}_{j,\ell}(\gamma)$  and $
U_n=\sqrt{n}\sum_{j=1}^{s}\sum_{\underset{\ell\neq j}{\ell=1}}^{s}\rho_\ell\left(\widehat{\Gamma}^{(n)}_{j,\ell}(\gamma)- \Gamma_{j,\ell}\right)$, with
\[
\widehat{\Gamma}^{(n)}_{j,\ell}(\gamma)=\parallel \widehat{m}_j \parallel_{\mathcal{H}}^2+\parallel  \widehat{m}_\ell\parallel_{\mathcal{H}}^2
-\frac{2}{n_j}\sum_{i=1}^{n_j}k_{i,n_j}(\gamma)<K(X_i^{(j)},\cdot),\widehat{m}_\ell>_{\mathcal{H}}
\]
and $\Gamma_{j,\ell}=\parallel  m_j-m_\ell\parallel_{\mathcal{H}}^2$.
Moreover
\begin{eqnarray*}
\delta_n
&=&\sqrt{n}\sum_{j=1}^{s}\sum_{\underset{\ell\neq j}{\ell=1}}^{s}(\pi_\ell-\rho_\ell)\,\left\{\left\|\widehat{m}_j\right\|^2_{\mathcal{H}}+\left\|\widehat{m}_\ell\right\|^2_{\mathcal{H}}-\frac{2}{n_j}\sum_{i=1}^{n_j}k_{i,n_j}(\gamma)<K(X_i^{(j)},\cdot),\widehat{m}_\ell>_{\mathcal{H}}\right\}\nonumber\\
&=&\sqrt{n}\sum_{j=1}^{s}\sum_{\underset{\ell\neq j}{\ell=1}}^{s}(\pi_\ell-\rho_\ell)\,\bigg\{\left\|\widehat{m}_j\right\|^2_{\mathcal{H}}+\left\|\widehat{m}_\ell\right\|^2_{\mathcal{H}}-\frac{2}{n_j}\sum_{i=1}^{n_j}\left(k_{i,n_j}(\gamma)-1\right)<K(X_i^{(j)},\cdot),\widehat{m}_\ell>_{\mathcal{H}}\nonumber\\
&&\hspace*{8cm}-\frac{2}{n_j}\sum_{i=1}^{n_j}<K(X_i^{(j)},\cdot),\widehat{m}_\ell>_{\mathcal{H}}\bigg\};
\end{eqnarray*}
since  we have
$
2\,n_j^{-1}\sum_{i=1}^{n_j}<K(X_i^{(j)},\cdot),\widehat{m}_\ell>_{\mathcal{H}}=2 <\widehat{m}_j,\widehat{m}_\ell>_{\mathcal{H}}$, putting $\widehat{\Gamma}^{n}_{j,\ell}=\parallel  \widehat{m}_j-\widehat{m}_\ell\parallel_{\mathcal{H}}^2$,  it follows
\begin{eqnarray*}\label{deltan}
\delta_n&=&\sqrt{n}\sum_{j=1}^{s}\sum_{\underset{\ell\neq j}{\ell=1}}^{s}(\pi_\ell-\rho_\ell)\,\bigg\{\widehat{\Gamma}^{n}_{j,\ell}-\frac{2}{n_j}\sum_{i=1}^{n_j}\left(k_{i,n_j}(\gamma)-1\right)<K(X_i^{(j)},\cdot),\widehat{m}_\ell>_{\mathcal{H}}\bigg\}\nonumber\\&=&A_n-2B_n,
\end{eqnarray*}
where $A_n$ and $B_n$ are the random variables given in (\ref{defan}) and (\ref{defbn}). Then, from Lemma \ref{lem1}, we deduce that $
\delta_n=o_P(1)$; thus $
\sqrt{n}\left(\widehat{\mathscr{T}}_{n,\gamma}-\mathscr{T}\right)=U_n+o_P(1)$. Therefore, it remains to get the asymptotic distribution of  $U_n$. We have
\begin{eqnarray}
U_n
&=&\sqrt{n}\sum_{j=1}^{s}\sum_{\underset{\ell\neq j}{\ell=1}}^{s}\rho_\ell\bigg\{\left\|\widehat{m}_j\right\|^2_{\mathcal{H}}+\left\|\widehat{m}_\ell\right\|^2_{\mathcal{H}}-\frac{2}{n_j}\sum_{i=1}^{n_j}\left(k_{i,n_j}(\gamma)-1\right)<K(X_i^{(j)},\cdot),\widehat{m}_\ell-m_\ell>_{\mathcal{H}}\nonumber\\
&&\hspace{2cm}-2<\widehat{m}_j,\widehat{m}_\ell>_{\mathcal{H}}+\frac{2}{n_j}\sum_{i=1}^{n_j}<K(X_i^{(j)},\cdot),m_\ell>_{\mathcal{H}}\nonumber\\
&&\hspace{2cm}-\frac{2}{n_j}\sum_{i=1}^{n_j}k_{i,n_j}(\gamma)<K(X_i^{(j)},\cdot),m_\ell>_{\mathcal{H}}- \left\|m_j-m_\ell\right\|^2_{\mathcal{H}}\bigg\}.\nonumber
\end{eqnarray}
Then, using the equalities
\begin{eqnarray*}
\left\|\widehat{m}_j\right\|^2_{\mathcal{H}}&=&\left\|\widehat{m}_j-m_j\right\|^2_{\mathcal{H}}+2<\widehat{m}_j,m_j>_{\mathcal{H}}-\left\|m_j\right\|^2_{\mathcal{H}}\\
&=&\left\|\widehat{m}_j-m_j\right\|^2_{\mathcal{H}}+\frac{2}{n_j}\sum_{i=1}^{n_j}<K(X_i^{(j)},\cdot),m_j>_{\mathcal{H}}-\left\|m_j\right\|^2_{\mathcal{H}}\nonumber\\
\end{eqnarray*}
and 
\begin{eqnarray}
<\widehat{m}_j,\widehat{m}_\ell>_{\mathcal{H}}&=& <\widehat{m}_j-m_j,\widehat{m}_\ell-m_\ell>_{\mathcal{H}}+\frac{1}{n_j}\sum_{i=1}^{n_j}<K(X_i^{(j)},\cdot),m_\ell>_{\mathcal{H}}\nonumber\\&&+\frac{1}{n_\ell}\sum_{i=1}^{n_\ell}<K(X_i^{(\ell)},\cdot),m_j>_{\mathcal{H}}-<m_j,m_\ell>_{\mathcal{H}},\nonumber
\end{eqnarray}
we obtain $U_n=-2C_n-2D_n+E_n+F_n$, where $C_n$ and $D_n$ are the random variables given in (\ref{defcn}) and (\ref{defdn}), 
$
E_n=\sqrt{n}\sum_{j=1}^{s}\sum_{\underset{\ell\neq j}{\ell=1}}^{s}\rho_\ell\left\{\left\|\widehat{m}_j-m_j\right\|^2_{\mathcal{H}}+\left\|\widehat{m}_\ell-m_\ell\right\|^2_{\mathcal{H}}\right\},
$
\begin{eqnarray}
&&F_n=\sqrt{n}\sum_{j=1}^{s}\sum_{\underset{\ell\neq j}{\ell=1}}^{s}\rho_\ell\bigg\{\frac{2}{n_j}\sum_{i=1}^{n_j}<K(X_i^{(j)},\cdot)-m_j,m_j>_{\mathcal{H}}\nonumber\\
&&\hspace*{3cm}+\frac{2}{n_\ell}\sum_{i=1}^{n_\ell}<K(X_i^{(\ell)},\cdot)-m_\ell,m_\ell>_{\mathcal{H}}\nonumber\\
&&\hspace*{3cm}-\frac{2}{n_\ell}\sum_{i=1}^{n_\ell}<K(X_i^{(\ell)},\cdot)-m_\ell,m_j>_{\mathcal{H}}\bigg\}\nonumber\\&&\hspace*{3cm}-\left.\frac{2}{n_j}\sum_{i=1}^{n_j}k_{i,n_j}(\gamma)<K(X_i^{(j)},\cdot)-m_j,m_\ell>_{\mathcal{H}}\right\}.\label{deffn}
\end{eqnarray}
From (\ref{diffmoy}) and the equality $\lim_{n_j\rightarrow +\infty}\frac{n_j}{n}=\rho_j$, we deduce that $E_n=o_P(1)$. This result  and  Lemma 1 imply that $U_n=F_n+o_P(1)$. Then, $U_n$ has the same limiting distribution than $F_n$ and it remains to derive this latter. Since   $\sum_{\underset{\ell\neq j}{\ell=1}}^{s}\rho_\ell=1-\rho_j$, we have
\begin{eqnarray}\label{f1}
\sum_{\underset{\ell\neq j}{\ell=1}}^{s}&\rho_\ell&\sum_{i=1}^{n_j}<K(X_i^{(j)},\cdot)-m_j,m_j>_{\mathcal{H}}\nonumber\\
&=&\sum_{i=1}^{n_j}<K(X_i^{(j)},\cdot)-m_j,(1-\rho_j)m_j>_{\mathcal{H}}.
\end{eqnarray}
Furthermore, 
\begin{eqnarray}\label{f2}
\sum_{j=1}^{s}\sum_{\underset{\ell\neq j}{\ell=1}}^{s}&\rho_\ell&\bigg\{\frac{2}{n_\ell}\sum_{i=1}^{n_\ell}<K(X_i^{(\ell)},\cdot)-m_\ell,m_\ell>_{\mathcal{H}}\bigg\}\nonumber\\
&=&\sum_{\ell=1}^{s}\sum_{\underset{j\neq \ell}{j=1}}^{s}\rho_\ell\bigg\{\frac{2}{n_\ell}\sum_{i=1}^{n_\ell}<K(X_i^{(\ell)},\cdot)-m_\ell,m_\ell>_{\mathcal{H}}\bigg\}\nonumber\\
&=&\sum_{\ell=1}^{s}\rho_\ell\bigg\{\frac{2}{n_\ell}\sum_{i=1}^{n_\ell}<K(X_i^{(\ell)},\cdot)-m_\ell,(s-1)m_\ell>_{\mathcal{H}}\bigg\}\nonumber\\
&=&\sum_{j=1}^{s}\frac{2}{n_j}\sum_{i=1}^{n_j}<K(X_i^{(j)},\cdot)-m_j,(s-1)\rho_j\,m_j>_{\mathcal{H}},
\end{eqnarray}
\begin{eqnarray}\label{f3}
\sum_{j=1}^{s}\sum_{\underset{\ell\neq j}{\ell=1}}^{s}&\rho_\ell&\bigg\{\frac{2}{n_\ell}\sum_{i=1}^{n_\ell}<K(X_i^{(\ell)},\cdot)-m_\ell,m_j>_{\mathcal{H}}\bigg\}\nonumber\\
&=&\sum_{\ell=1}^{s}\sum_{\underset{j\neq \ell}{j=1}}^{s}\rho_\ell\bigg\{\frac{2}{n_\ell}\sum_{i=1}^{n_\ell}<K(X_i^{(\ell)},\cdot)-m_\ell,m_j>_{\mathcal{H}}\bigg\}\nonumber\\
&=&\sum_{\ell=1}^{s}\frac{2}{n_\ell}\sum_{i=1}^{n_\ell}<K(X_i^{(\ell)},\cdot)-m_\ell,\rho_\ell (\mu-m_\ell)>_{\mathcal{H}}\nonumber\\
&=&\sum_{j=1}^{s}\frac{2}{n_j}\sum_{i=1}^{n_j}<K(X_i^{(j)},\cdot)-m_j,\rho_j (\mu-m_j)>_{\mathcal{H}}
\end{eqnarray}
and, since  $\sum_{\underset{\ell\neq j}{\ell=1}}^{s}\rho_\ell\,m_\ell=m-\rho_j\,m_j$, 
\begin{eqnarray}\label{f4}
 \sum_{\underset{\ell\neq j}{\ell=1}}^{s}&\rho_\ell&\bigg\{\frac{2}{n_j}\sum_{i=1}^{n_j}k_{i,n_j}(\gamma)<K(X_i^{(j)},\cdot)-m_j,m_\ell>_{\mathcal{H}}\bigg\}\nonumber\\
&=&\frac{2}{n_j}\sum_{i=1}^{n_j}k_{i,n_j}(\gamma)<K(X_i^{(j)},\cdot)-m_j,m-\rho_j\,m_j>_{\mathcal{H}}.
\end{eqnarray}
Then, using (\ref{deffn}), (\ref{f1}),  (\ref{f2}),  (\ref{f3}),  (\ref{f4}) and the equality $\sqrt{n}\,n_j^{-1}=\pi_j^{-1/2}\,n_j^{-1/2}$, we obtain
\begin{equation}\label{fn}
F_n=\sum_{j=1}^{s}2\pi_j^{-1/2}Y_{n,j,\gamma} ,
\end{equation} 
where
$
Y_{n,j,\gamma}=\frac{1}{\sqrt{n_j}}\sum_{i=1}^{n_j}\,\bigg\{\mathcal{U}_{j}(X_i^{(j)})-k_{i,n_j}(\gamma)\mathcal{V}_{j}(X_i^{(j)})\bigg\}.
$
Let us put  $s_{n,j,\gamma}^2=\sum_{i=1}^{n_j}Var\left(\mathcal{W}_{n,i,j}(X_i^{(j)})\right)$, where $\mathcal{W}_{n,i,j}(X_i^{(j)})=\mathcal{U}_{j}(X_i^{(j)})-k_{i,n_j}(\gamma)\mathcal{V}_{j}(X_i^{(j)})$. By similar arguments than in the proof of Theorem 1 in Makigusa and Naito (2020) we obtain that, for any $\varepsilon>0$,
\[
s_{n,j,\gamma}^{-2}\sum_{i=1}^{n_j}\int_{\{x:| \mathcal{U}_{j}(x)-k_{i,n_j}(\gamma)\mathcal{V}_{j}(x) |>\varepsilon s_{n,j,\gamma}\}}^{} \bigg(\mathcal{U}_{j}(x)-k_{i,n_j}(\gamma)\mathcal{V}_{j}(x)\bigg)^2\,d\mathbb{P}_{j}(x) 
\]
converges to $0$ as $n_j\rightarrow +\infty$. Therefore, by Section 1.9.3 in Serfling (1980) we obtain that  
 $s_{n,j,\gamma}^{-1}\sum_{i=1}^{n_j}\mathcal{W}_{n,i,j}(X_i^{(j)})	\stackrel{\mathscr{D}}{\rightarrow} \mathcal{N}\left(0,1\right)$, that is 
$\sqrt{n_j}s_{n,j,\gamma}^{-1}Y_{n,j,\gamma}	\stackrel{\mathscr{D}}{\rightarrow} \mathcal{N}\left(0,1\right)$. However,
\begin{eqnarray*}
\left(\frac{s_{n,j,\gamma}}{\sqrt{n_j}}\right)^2&=&\frac{1}{n_j}\sum_{i=1}^{n_j}\bigg\{Var\left(\mathcal{U}_{j}(X_i^{(j)})\right)+k^2_{i,n_j}(\gamma)Var\left(\mathcal{V}_{j}(X_i^{(j)})\right)\nonumber\\
&&-2k_{i,n_j}(\gamma)\,Cov\left(\mathcal{U}_{j}(X_i^{(j)}),\mathcal{V}_{j}(X_i^{(j)})\right)\bigg\}\nonumber\\
&=&Var\left(\mathcal{U}_{j}(X_1^{(j)})\right)+\left(\frac{1}{n_j}\sum_{i=1}^{n_j}k^2_{i,n_j}(\gamma)\right)Var\left(\mathcal{V}_{j}(X_1^{(j)})\right)\nonumber\\
&&-2\left(\frac{1}{n_j}\sum_{i=1}^{n_j}k_{i,n_j}(\gamma)\right)Cov\left(\mathcal{U}_{j}(X_1^{(j)}),\mathcal{V}_{j}(X_1^{(j)})\right)\bigg\};
\end{eqnarray*}
then, using   $ (\mathscr{A}_3) $ and $ (\mathscr{A}_5) $, we get $\lim\limits_{n_j\rightarrow +\infty}\left(n_j^{-1}s_{n,j,\gamma}^2\right)=\sigma_j^2(\gamma)$, where 
$\sigma_j^2(\gamma)$ is defined in (\ref{sigmaj}). Hence, $Y_{n,j,\gamma}	\stackrel{\mathscr{D}}{\rightarrow} \mathcal{N}\left(0,\sigma_j^2(\gamma)\right)$. Since  $Y_{n,j,\gamma}$ and $Y_{n,\ell,\gamma}$ are independent when $j\neq\ell$, we deduce from (\ref{fn}) and the equality  $\lim\limits_{n_j\rightarrow +\infty}(\pi_j)=\rho_j$ that  $F_n\stackrel{\mathscr{D}}{\rightarrow} \mathcal{N}\left(0,\sigma_\gamma^2\right)$, where  $\sigma^2_{\gamma}=\sum_{j=1}^{s}4\rho_j^{-1}\sigma_j^2(\gamma)$.
\subsection{Proof of Proposition \ref{pro}}
\noindent It suffices to prove that $\widehat{\nu}_j^2$ is consistent for $\nu^2$. First, by Cauchy-Schwartz inequality and (\ref{kinf}),
\begin{eqnarray*}
\frac{1}{n_j}\sum_{i=1}^{n_j}<K(X_i^{(j)},\cdot) ,\widehat{m}-m>_{\mathcal{H}}^2
&\leq&\frac{1}{n_j}\sum_{i=1}^{n_j}\Vert K(X_i^{(j)},\cdot) \Vert_{\mathcal{H}}^2\,\Vert\widehat{m}-m\Vert_{\mathcal{H}}^2\\
&\leq&\frac{1}{n}\Vert K \Vert_{\infty}\,\Vert\sqrt{n}(\widehat{m}-m)\Vert_{\mathcal{H}}^2.
\end{eqnarray*}
Using (\ref{diffmoy}) and ($\mathscr{A}_2$) it is easy to check that $\Vert\sqrt{n}(\widehat{m}-m)\Vert_{\mathcal{H}}=O_P(1)$. Hence, from the previous inequality, $n_j^{-1}\sum_{i=1}^{n_j}<K(X_i^{(j)},\cdot) ,\widehat{m}-m>_{\mathcal{H}}^2=o_P(1)$. Another use of Cauchy-Schwartz inequality and (\ref{kinf}) gives the inequality
\begin{eqnarray*}
\bigg\vert\frac{1}{n_j}\sum_{i=1}^{n_j}<K(X_i^{(j)},\cdot) ,\widehat{m}-m>_{\mathcal{H}}<K(X_i^{(j)},\cdot) ,m>_{\mathcal{H}}\bigg\vert
&\leq&\frac{\Vert m \Vert_{\mathcal{H}}\Vert K \Vert_{\infty}}{\sqrt{n}}\,\Vert\sqrt{n}(\widehat{m}-m)\Vert_{\mathcal{H}}
\end{eqnarray*}
that implies $n_j^{-1}\sum_{i=1}^{n_j}<K(X_i^{(j)},\cdot) ,\widehat{m}-m>_{\mathcal{H}}<K(X_i^{(j)},\cdot) ,m>_{\mathcal{H}}=o_P(1)$. Thus, $n_j^{-1}\sum_{i=1}^{n_j}<K(X_i^{(j)},\cdot) ,\widehat{m}>_{\mathcal{H}}^2=o_P(1)+n_j^{-1}\sum_{i=1}^{n_j}<K(X_i^{(j)},\cdot) ,m>_{\mathcal{H}}^2$, and using the law of large numbers we conclude that $n_j^{-1}\sum_{i=1}^{n_j}<K(X_i^{(j)},\cdot) ,\widehat{m}>_{\mathcal{H}}^2$ converges in probability, as $n_j\rightarrow +\infty$,  to $\mathbb{E}\left(<K(X_1^{(j)},\cdot) ,m>_{\mathcal{H}}^2\right)$. On the other hand, Cauchy-Schwartz inequality and (\ref{kinf}) give the inequality
\begin{eqnarray*}
\bigg\vert\frac{1}{n_j}\sum_{i=1}^{n_j}<K(X_i^{(j)},\cdot) ,\widehat{m}-m>_{\mathcal{H}}\bigg\vert
&\leq&\frac{ \Vert K \Vert_{\infty}^{1/2}}{\sqrt{n}}\,\Vert\sqrt{n}(\widehat{m}-m)\Vert_{\mathcal{H}}
\end{eqnarray*}
that implies $n_j^{-1}\sum_{i=1}^{n_j}<K(X_i^{(j)},\cdot) ,\widehat{m}-m>_{\mathcal{H}}=o_P(1)$. Therefore, we have  $n_j^{-1}\sum_{i=1}^{n_j}<K(X_i^{(j)},\cdot) ,\widehat{m}>_{\mathcal{H}}=o_P(1)+n_j^{-1}\sum_{i=1}^{n_j}<K(X_i^{(j)},\cdot) ,m>_{\mathcal{H}}$; from  the law of large numbers we deduce that $n_j^{-1}\sum_{i=1}^{n_j}<K(X_i^{(j)},\cdot) ,\widehat{m}>_{\mathcal{H}}$ converges in probability, as $n_j\rightarrow +\infty$,  to $\mathbb{E}\left(<K(X_1^{(j)},\cdot) ,m>_{\mathcal{H}}\right)$. The preceding convergences properties imply that  $\widehat{\nu}_j^2$ converges in probability, as $n_j\rightarrow +\infty$,  to $\nu_j^2:=\mathbb{E}\left(<K(X_1^{(j)},\cdot) ,m>_{\mathcal{H}}^2\right)-\mathbb{E}\left(<K(X_1^{(j)},\cdot) ,m>_{\mathcal{H}}\right)^2$. Under the hypothesis  $\mathbb{P}_1=\cdots=\mathbb{P}_s$, we have $\nu_j^2=\nu^2$.

\end{document}